\documentclass[12pt,a4paper]{amsart}
\usepackage{enumerate}
\usepackage{amssymb}
\usepackage{amsmath}
\usepackage{mathrsfs}
\usepackage{amsthm}
\usepackage{mathtools}
\theoremstyle{definition}
\newtheorem{theorem}{Theorem}
\newtheorem{lemma}[theorem]{Lemma}

\newtheorem{remark}[theorem]{Remark}
\usepackage[colorlinks,linkcolor=blue]{hyperref}
\begin{document}
\title[Level sets of harmonic functions]{Level sets of harmonic functions on three-dimensional manifolds with nonnegative scalar curvature}
\author{Yukai Sun}
\address{School of Mathematics and Statistics, Henan University, Kaifeng 475004 P. R. China and Center for Applied Mathematics of Henan Province, Henan University, Zhengzhou 450046 P. R. China}
\email{sunyukai@henu.edu.cn}
\date{}
\begin{abstract}
  We investigate the level sets of harmonic functions on $(\mathbb{R}^{3}\setminus \{0\},g)$. Drawing inspiration from Miao\cite{M-SCM-2025}, we adopt the method developed by Munteanu-Wang\cite{WM-IMRN-2023} to derive a monotonic quantity associated with the level sets of harmonic functions on $(\mathbb{R}^{3}\setminus \{0\},g)$ with nonnegative scalar curvature, under certain conditions. Furthermore, we establish a rigidity result for this quantity. Additionally, we find an extra scalar-flat metric on $\mathbb{R}^{3}\setminus \{0\}$.
\end{abstract}
\maketitle
\section{Introduction}
In recent years, extensive research has been carried out on the level sets of harmonic functions over three-dimensional manifolds. Miao \cite{M-PekingMath-2025, M-SCM-2025} derived monotonicity formulas for asymptotically flat manifolds with nonnegative scalar curvature, drawing on the analysis of harmonic function level sets. Notably, these monotonicity formulas are capable of characterizing both Euclidean space and the Schwarzschild metric. Additionally, Munteanu-Wang\cite{WM-IMRN-2023} established a monotonicity formula for three-dimensional manifolds with nonnegative scalar curvature and provided a characterization of its rigidity case. For more relevant results, refer to Agostiniani-Mazzieri-Mazzieri-Oronzio\cite{AMMO-ANL-2023}, Agostiniani-Mazzieri-Oronzio\cite{AMO-CMP-2024}, Bray-Kazaras-Khuri-Stern\cite{BKKS-JGA-2022}, Colding\cite{Colding-Acta-2012}, Colding-Minicozzi\cite{CM-PNAS-2013,CM-CVPDE-2014}, Munteanu-Wang\cite{MW-JFA-2024}, Stern \cite{Stern-JDG-2022}.

Motivated by the work of Miao\cite{M-SCM-2025} and Munteanu-Wang\cite{WM-IMRN-2023}, we investigate three-manifolds $M^{3}$ that are diffeomorphic to $(0,\infty)\times \mathbb{S}^2$.

Let $(M^{3},g)$ be an oriented three-dimensional manifold with nonnegative scalar curvature $\operatorname{Sc}_{g}$. We assume $M=\mathbb{R}^{3}\setminus\{0\}$ . Let $G$ be a harmonic function on $(M,g)$, its level sets $l(t)$ are defined as
\[l(t)=\{x\in M| G(x)=t\}.\]
We consider the following quantity
\[w(t)=\int_{l(t)}|\nabla G|^2.\]
Hereafter, the area element is omitted for brevity.

\begin{theorem}\label{Thm-c=1}
  Let $(M,g)$ as above and $f(t)=t^{-2}$. Assume $G$ is a harmonic function on $(M^{3},g)$ with the following properties:
  \begin{enumerate}
    \item $G(x)\to 1$ as $x\to 0$ and $G(x)\to 0$ as $x\to \infty$;
    \item \(\lim_{t \to 1} \int_{l(t)} f(G) \frac{\left\langle \nabla |\nabla G|, \nabla G \right\rangle}{|\nabla G|} = 0\) and

     \(\lim_{t \to 1} \int_{l(t)} |\nabla G| \frac{\left\langle \nabla f(G), \nabla G \right\rangle}{|\nabla G|} = 0\);
    \item $\lim_{t\to 0}\frac{w(t)}{t}=0$ and $\lim_{t\to 1}\frac{w(t)}{1-t}=0$.
    \item For each $t$, $l(t)$ is connected.
  \end{enumerate}
  Then, for $t\in (0,1)$,
  \begin{equation}\label{eqn-similar-to-miao}
    \frac{d\left[(1-t)^{-3}t^{-1}w(t)-4\pi (1-t)^{-1}\right]}{dt}\leq 0.
  \end{equation}
   If
  \[\frac{d\left[(1-t)^{-3}t^{-1}w(t)-4\pi (1-t)^{-1}\right]}{dt}=0\]
  for all $t\in (0,1)$, then
  \[g=\frac{c^2}{t^4(1-t)^4}dt^2+\frac{c^2}{t^2(1-t)^2}g_{\mathbb{S}^2},\]
  for some constant $c>0$, where $g_{\mathbb{S}^2}$ is the metric of sectional curvature $1$ on $\mathbb{S}^2$.
\end{theorem}
\begin{remark}
  The monotonicity formula in equation (\ref{eqn-similar-to-miao}) is similar to Theorem 1.3 in \cite{M-SCM-2025} with different conditions. This metric $g$ is the Schwarzschild metric $\left(1+\frac{c}{|x|}\right)^4g_{\mathbb{R}^{3}}$ on $\mathbb{R}^{3}\setminus \{0\}$. To see that you can set
  \[t = 1-\left(1+\frac{c}{r}\right)^{-1}\]
  for $r=|x|>0$. Then, through a direct computation, one can obtain the Schwarzschild metric on $\mathbb{R}^{3}\setminus \{0\}$. The function
  \[G(x)=1-\left(1+\frac{c}{|x|}\right)^{-1}\]
  is harmonic with respect to the metric $g$ and satisfies the above conditions.
\end{remark}
We adopt the method developed by Munteanu-Wang \cite{WM-IMRN-2023} to prove Theorem \ref{Thm-c=1}, though there exist discrepancies in certain details. Notably, the approaches employed to establish the rigid part of the theorem are entirely distinct. Moreover, our proof suggests that Munteanu-Wang’s method may be extendable to construct a variety of scalar-flat metrics; further discussion on this can be found in Remark \ref{Rmk-fh}.

The result of Theorem \ref{Thm-c=1} is analogous to that of Munteanu-Wang\cite{WM-IMRN-2023}, where their result characterizes Euclidean space. In contrast, our result serves to characterize the Schwarzschild metric. Here is their result.
\begin{theorem}[Theorem 1.1 \cite{WM-IMRN-2023}]
  Let $(M, g)$ be a complete noncompact 3D manifold with nonnegative scalar curvature. Assume that $M$ has one end and its 1st Betti number $b_1(M) = 0$. If $M$ is nonparabolic and the minimal positive Green's function $G(x) = G(p, x)$ satisfies $\lim_{x \to \infty} G(x) = 0$, then
\[
\frac{\mathrm{d}}{\mathrm{d}t} \left( \frac{1}{t} \int_{l(t)} |\nabla G|^2 - 4\pi t \right) \leq 0
\]
for all regular values $t > 0$. Moreover, equality holds for some $T > 0$ if and only if the super level set $\{ x \in M : G(x) > T \}$ is isometric to a ball in the Euclidean space $\mathbb{R}^3$.
\end{theorem}

If we change the initial condition about the harmonic function $G$, there is another result.

\begin{theorem}\label{Thm-c=-1}
  Let $(M,g)$ as before and $f(t)=t^{-2}$. Assume $G$ is a harmonic function on $(M^{3},g)$ with the following properties:
  \begin{enumerate}
    \item $G(x)\to -1$ as $x\to 0$ and $G(x)\to -\infty$ as $x\to \infty$;
    \item \(\lim_{t \to -1} \int_{l(t)} f(G) \frac{\left\langle \nabla |\nabla G|, \nabla G \right\rangle}{|\nabla G|} = 0\) and

     \(\lim_{t \to -1} \int_{l(t)} |\nabla G| \frac{\left\langle \nabla f(G), \nabla G \right\rangle}{|\nabla G|} = 0\);
    \item $w(2)=16\pi$ and $\lim_{t\to -1}\frac{w(t)}{1+t}=0$.
    \item For each $t$, $l(t)$ is connected.
  \end{enumerate}
  Then, for $t\in (-\infty,-1)$,
  \[\frac{d\left[(1+t)^{-3}t^{-1}w(t)+4\pi (1+t)^{-1}\right]}{dt}\leq 0.\]
   If
  \[\frac{d\left[(1+t)^{-3}t^{-1}w(t)+4\pi (1+t)^{-1}\right]}{dt}=0\]
  for all $t\in (-\infty,-1)$, then
  \[g=\frac{c}{t^4(1+t)^4}dt^2+\frac{c}{t^2(1+t)^2}g_{\mathbb{S}^2},\text{ on } (-\infty,-1)\times\mathbb{S}^2, \text{ and } c>0,\]
  where $g_{\mathbb{S}^2}$ is the metric of sectional curvature $1$ on $\mathbb{S}^2$. The metric $g$ is a scalar flat metric.
\end{theorem}
\begin{remark}
    Let $r=|x|=t(1+t)$ on $\mathbb{R}^{3}\setminus\{0\}$, then the metric $g$ in Theorem \ref{Thm-c=-1} becomes
    \[g=cr^{-4}(1+4r)^{-1}\left[dr^2+r^2(1+4r)g_{\mathbb{S}^{2}}\right], \text{ on } (0,\infty)\times\mathbb{S}^2.\]
    The function
    \[G(x)=-\frac{1}{2}-\frac{1}{2}\sqrt{1+4|x|}\]
    is harmonic with respect to the metric $g$ and satisfies the properties as in Theorem \ref{Thm-c=-1}.
\end{remark}

The paper is organized as follows. In section \ref{sec-basic}, we give some useful Lemmas. In section \ref{sec-proof}, we proof Theorem \ref{Thm-c=1} and \ref{Thm-c=-1}.

\section{Basic facts}\label{sec-basic}
In this section, we present several key facts concerning the level sets of harmonic functions. These results are taken from the work of Munteanu-Wang \cite{WM-IMRN-2023}.
\begin{lemma}[Lemma 2.1 \cite{WM-IMRN-2023}]\label{lem-Ric}
  Let $(M^{3},g)$ be a 3D complete noncompact Riemannian manifold with scalar curvature $\operatorname{Sc}_{g}$ and $u$ a harmonic function on $M$. Then on each regular level set $l(r)$ of $u$,
\[
\operatorname{Ric}_{g}(\nabla u, \nabla u) |\nabla u|^{-2} = \frac{1}{2}\operatorname{Sc}_{g} - \frac{1}{2}\operatorname{Sc}_{l(r)} + \frac{1}{|\nabla u|^2} \left( \left| \nabla |\nabla u| \right|^2 - \frac{1}{2} \left| \nabla^2 u \right|^2 \right),
\]
where $\operatorname{Sc}_{l(r)}$ denotes the scalar curvature of $l(r)$.
\end{lemma}

On a regular level set $l(r)$ of $u$, its unit normal vector is given by
\[e_{1}=\frac{\nabla u}{|\nabla u|}.\]
Choose $\{e_{a}\}_{a=2,3}$, unit vectors tangent to $l(r)$, such that $\{e_{1},e_{2},e_{3}\}$ forms a local orthonormal frame on $M$. Since $u$ is harmonic, the second fundamental form $A_{r}$ and the mean curvature $H_{r}$ of $l(r)$ are given by
\[(A_{r})_{ab}=\frac{u_{ab}}{|\nabla u|} \text{ and } H_{r}=-\frac{u_{11}}{|\nabla u|}, \text{ respectively }.\]

\begin{lemma}[Lemma 2.4\cite{WM-IMRN-2023}]\label{lem-ineq-u}
  Let $(M,g)$ be a 3-dimensional complete noncompact Riemannian manifold and $u$ a harmonic function on $M$. Then
  \[|\nabla^2u|^2\geq \frac{3}{2}|\nabla |\nabla u||^2 \text{ on } M.\]
\end{lemma}
We write down the proof to analyze the equality case.
\begin{proof}
  It suffices to prove this at points where $|\nabla u|\neq 0$. Let \[e_{1}=\frac{\nabla u}{|\nabla u|}.\]
Choose $\{e_{a}\}_{a=2,3}$ such that $\{e_{1},e_{2},e_{3}\}$ is a local orthonormal frame on $M$. Then
\begin{eqnarray*}
  |u_{ij}|^2 &=& |u_{11}|^2+\sum_{a=1}^{2}2|u_{1a}|^2+\sum_{a,b=1}^{2}|u_{ab}|^2\\
  &\geq&|u_{11}|^2+\sum_{a=1}^{2}2|u_{1a}|^2+\frac{1}{2}|u_{22}+u_{33}|^2\\
  &=&\frac{3}{2}|u_{11}|^2+\sum_{a=1}^{2}2|u_{1a}|^2,
\end{eqnarray*}
since $u$ is harmonic. Therefore
\begin{eqnarray*}
  |u_{ij}|^2 &\geq& \frac{3}{2}\left(|u_{11}|^2+\sum_{a=1}^{2}|u_{1a}|^2\right)\\
  &=&\frac{3}{2}|\nabla|\nabla u||^2.
\end{eqnarray*}
\end{proof}

\section{Proof of Theorems}\label{sec-proof}
We employ the method developed by Munteanu-Wang \cite{WM-IMRN-2023} to prove Theorem \ref{Thm-c=1}. For the sake of a complete analysis, we elaborate on the key details herein.

Recall that $l(t) = \{x\in M |G(x)=t\}$. Let $L(a,b) = \{x\in M |a<G(x)<b\}.$ Consider the function
\[w(t)=\int_{l(t)}|\nabla G|^2.\]
Whenever $l(t)$ is regular, its mean curvature $H_{l(t)}$ is given by
\[H_{l(t)}=\frac{\sum_{a}G_{aa}}{|\nabla G|}=-\frac{G_{11}}{|\nabla G|}=-\frac{\langle\nabla |\nabla G|,\nabla G\rangle}{|\nabla G|^2},\]
where $e_{1}=\frac{\nabla G}{|\nabla G|}$ and $\{e_{a}\}_{a=2,3}$ are unit tangent vectors on $l(t)$ such that $\{e_{1},e_{2},e_{3}\}$ is a local orthonormal frame on $M$. We compute
\begin{eqnarray}\label{eqn-dw-w}
  \frac{d w(t)}{dt}& = &\int_{l(t)}\left(\frac{\langle\nabla|\nabla G|^2,\nabla G\rangle}{|\nabla G|^2}+\frac{H_{l(t)}}{|\nabla G|}|\nabla G|^2\right)\\\nonumber
  &=&\int_{l(t)}\left(\frac{\langle\nabla|\nabla G|,\nabla G\rangle}{|\nabla G|}\right).
\end{eqnarray}
Multiplying the equation by $f(t)$, we obtain
\[f(t)\frac{d w(t)}{dt}=\int_{l(t)}\left(\frac{\langle\nabla|\nabla G|,\nabla G\rangle}{|\nabla G|}\right)f(G)\]
On the other hand, by Green’s identity, we have
\begin{eqnarray*}
   && \int_{L(t,T)}\left(f(G)\Delta |\nabla G|-|\nabla G|\Delta f(G)\right)\\
   &=&\int_{l(T)}\left(f(G)\frac{\langle\nabla|\nabla G|,\nabla G\rangle}{|\nabla G|}-|\nabla G|\frac{\langle\nabla f(G),\nabla G\rangle}{|\nabla G|}\right)\\
   &&-\int_{l(t)}\left(f(G)\frac{\langle\nabla|\nabla G|,\nabla G\rangle}{|\nabla G|}-|\nabla G|\frac{\langle\nabla f(G),\nabla G\rangle}{|\nabla G|}\right).
\end{eqnarray*}
Since $G$ is harmonic on $L(t,T)$, we have
\[|\nabla G|\Delta f(G)=f''(G)|\nabla G|^{3}.\]

Since
\[
\lim_{T \to 1} \int_{l(T)} f(G) \frac{\left\langle \nabla |\nabla G|, \nabla G \right\rangle}{|\nabla G|} = 0\]
and
\[ \lim_{T \to 1} \int_{l(T)} |\nabla G| \frac{\left\langle \nabla f(G), \nabla G \right\rangle}{|\nabla G|} = 0.
\]

Then,
\[
\left| \int_{L(t,1)} \left( f(G) \Delta |\nabla G| - |\nabla G| \Delta f(G) \right) \right| < \infty,
\]
and the following identity holds
\begin{eqnarray*}
   && \int_{l(t)}\left(f(G)\frac{\langle\nabla|\nabla G|,\nabla G\rangle}{|\nabla G|}-|\nabla G|\frac{\langle\nabla f(G),\nabla G\rangle}{|\nabla G|}\right)\\
   &=&-\int_{L(t,1)}\left(f(G)\Delta |\nabla G|-|\nabla G|\Delta f(G)\right).
\end{eqnarray*}

We conclude that

\begin{eqnarray*}
  f(t)\frac{d w(t)}{dt}&=&\int_{l(t)}\left(\frac{\langle\nabla|\nabla G|,\nabla G\rangle}{|\nabla G|}\right)f(G)\\
   &=& \int_{l(t)}f'(G)|\nabla G|^2-\int_{L(t,1)}f(G)\Delta |\nabla G|\\
   &&+\int_{L(t,1)}f''(G)|\nabla G|^{3}.
\end{eqnarray*}

Note that by the co-area formula
\begin{eqnarray*}
  \int_{L(t,1)}f''(G)|\nabla G|^{3} &=& \int_{t}^{1}f''(t)\int_{l(r)}|\nabla G|^{2}\\
  &=&\int_{t}^{1}f''(r)w(r)dr.
\end{eqnarray*}

Hence, we have
\begin{eqnarray*}
  f(t)\frac{d w(t)}{dt}= f'(t)w(t)+\int_{t}^{1}f''(r)w(r)dr-\int_{L(t,1)}f(G)\Delta |\nabla G|.
\end{eqnarray*}

We now estimate the last term. Using the Bochner formula
\[  \Delta |\nabla G| = \left(|\operatorname{Hess}G|^2-|\nabla |\nabla G||^2\right)|\nabla G|^{-1}+\operatorname{Ric}_{g}(\nabla G,\nabla G)|\nabla G|^{-1}.\]

We have
\begin{eqnarray*}
  \int_{l(r)}|\nabla G|^{-1}\Delta |\nabla G| &=& \int_{l(r)}\left(|\operatorname{Hess}G|^2-|\nabla |\nabla G||^2\right)|\nabla G|^{-2}\\
  &&+\int_{l(r)}\operatorname{Ric}_{g}(\nabla G,\nabla G)|\nabla G|^{-2}.
\end{eqnarray*}

Applying Lemma \ref{lem-Ric} to $G$ gives
\begin{eqnarray*}
  \operatorname{Ric}_{g}(\nabla G, \nabla G) |\nabla G|^{-2} &=& \frac{1}{2}\operatorname{Sc}_{g} - \frac{1}{2}\operatorname{Sc}_{l(r)}\\
   &&+ \frac{1}{|\nabla G|^2} \left( \left| \nabla |\nabla G| \right|^2 - \frac{1}{2} \left| \nabla^2 G \right|^2 \right).
\end{eqnarray*}
We therefore conclude that
\begin{eqnarray*}
  \int_{l(r)}|\nabla G|^{-1}\Delta |\nabla G| &=& \frac{1}{2}\int_{l(r)}\left(|\operatorname{Hess}G|^2|\nabla G|^{-2}+\operatorname{Sc}_{g} - \operatorname{Sc}_{l(r)}\right).
\end{eqnarray*}
By Lemma \ref{lem-ineq-u},
\[|\operatorname{Hess}G|^2\geq \frac{3}{2}|\nabla|\nabla G||^2.\]
Since $l(r)$ is compact and connected for any $1>r>0$, the Gauss-Bonnet theorem implies that
\[\int_{l(r)}\operatorname{Sc}_{l(r)}=4\pi\chi(l(r))\leq 8\pi\]
whenever $r$ is a regular value of $G$.  Therefore, on any regular level set $l(r)$, one obtains that

\begin{eqnarray}\label{eqn-scalar}
  \int_{l(r)}|\nabla G|^{-1}\Delta |\nabla G| &\geq& \frac{3}{4}\int_{l(r)}|\nabla|\nabla G||^2|\nabla G|^{-2}-4\pi.
\end{eqnarray}

From equation (\ref{eqn-dw-w}), by the Cauchy-Schwarz inequality, we have
\begin{eqnarray}\label{eqn-w'}
  \left(w'(r)\right)^2 &=& \left(\int_{l(r)}\frac{\langle\nabla|\nabla G|,\nabla G\rangle}{|\nabla G|}\right)^2\\\nonumber
  &\leq&\left(\int_{l(r)}|\nabla|\nabla G||^2|\nabla G|^{-2}\right)\left(\int_{l(r)}|\nabla G|^2\right),
\end{eqnarray}

which says that
\[ \left(\int_{l(r)}|\nabla|\nabla G||^2|\nabla G|^{-2}\right)\geq (w'(r))^2w^{-1}(r) .\]

We conclude that
\begin{eqnarray*}
  \int_{l(r)}|\nabla G|^{-1}\Delta |\nabla G|\geq\frac{3}{4}(w'(r))^2w^{-1}(r)-4\pi.
\end{eqnarray*}

By the co-area formula, it follows that
\begin{eqnarray*}
 -\int_{L(t,1)}f(G)\Delta |\nabla G|&=& -\int_{t}^{1}f(r)\int_{l(r)}|\nabla G|^{-1}\Delta|\nabla G|\\
 &\leq&-\frac{3}{4}\int_{t}^{1}f(r)(w'(r))^2w^{-1}(r)dr\\
 &&+4\pi\int_{t}^{1}f(r).
\end{eqnarray*}

Then,
\begin{eqnarray*}
  f(t)w'(t)&=& f'(t)w(t)+\int_{t}^{1}f''(r)w(r)dr-\int_{L(t,1)}f(G)\Delta |\nabla G|\\
  &\leq& f'(t)w(t)+\int_{t}^{1}f''(r)w(r)dr\\
  &&-\frac{3}{4}\int_{t}^{1}f(r)(w'(r))^2w^{-1}(r)dr+4\pi\int_{t}^{1}f(r)
\end{eqnarray*}

Let $h(r)=\frac{1-2r}{r(1-r)}$. From the elementary inequality
\begin{eqnarray}\label{eqn-h}
  0 &\leq& w\left(\frac{w'}{w}-2h\right)^2\\\nonumber
  &=&\frac{(w')^2}{w}-4hw'+4h^2w.
\end{eqnarray}
One sees that
\begin{eqnarray*}
&&-\frac{3}{4}\int_{t}^{1}f(r)(w'(r))^2w^{-1}(r)dr\\
&\leq &-3\int_{t}^{1}f(r)h(r)w'(r)dr+3\int_{t}^{1}f(r)h^2(r)w(r)dr\\
&=&-3f(r)h(r)w(r)|_{t}^{1}+3\int_{t}^{1}(f'(r)h(r)+f(r)h'(r))w(r)dr\\
&&+3\int_{t}^{1}f(r)h^2(r)w(r)dr.
\end{eqnarray*}

Then,
\begin{eqnarray*}
  f(t)w'(t)&\leq& f'(t)w(t)+\int_{t}^{1}f''(r)w(r)dr\\
  &&-\frac{3}{4}\int_{t}^{1}f(r)(w'(r))^2w^{-1}(r)dr+4\pi\int_{t}^{1}f(r)\\
  &\leq&f'(t)w(t)+\int_{t}^{1}f''(r)w(r)dr\\
  &&+4\pi\int_{t}^{1}f(r)dr+3f(t)h(t)w(t)\\
  &&+3\int_{t}^{1}(f'(r)h(r)+f(r)h'(r))w(r)dr\\
  &&+3\int_{t}^{1}f(r)h^2(r)w(r)dr.
\end{eqnarray*}

 Recall that $f(r)=r^{-2}$. A direct computation obtains
\begin{equation}\label{eqn-f-h}
  f''(r)+3(f'(r)h(r)+f(r)h'(r))+3f(r)h^2(r)=0.
\end{equation}

Then
\begin{eqnarray*}
  t^{-2} w'(t)&\leq& -2t^{-3} w(t)+4\pi(t^{-1}-1)+3t^{-2}\frac{1-2t}{t(1-t)}w(t)\\
  &=&\left(-2+3\frac{1-2t}{1-t}\right)t^{-3} w(t)+4\pi(t^{-1}-1)\\
  &=&\left(\frac{1-4t}{1-t}\right)t^{-3} w(t)+4\pi(t^{-1}-1).
\end{eqnarray*}
And
\begin{eqnarray*}
  t(1-t) w'(t)&\leq&\left(1-4t\right) w(t)+4\pi(1-t)^2t^{2}\\
 \left( \frac{1-t}{t}w(t)\right)'&\leq&4\pi(1-t)^2-\frac{4}{(1-t)}\frac{(1-t)w(t)}{t}.
\end{eqnarray*}

Hence
\begin{equation*}
  \left[(1-t)^{-4}\left( \frac{1-t}{t}w(t)\right)-4\pi(1-t)^{-1}\right]'\leq 0.
\end{equation*}
I.e.,
\begin{equation*}
  \left[(1-t)^{-3}t^{-1}w(t)-4\pi(1-t)^{-1}\right]'\leq 0.
\end{equation*}

If
\[\left[(1-t)^{-4}\left( \frac{1-t}{t}w(t)\right)-4\pi(1-t)^{-1}\right]'=0,\;\forall t\in (0,1).\]

By the condition,
\[\lim_{t\to 0}\frac{w(t)}{t}=0.\]
We have
\[\left[(1-t)^{-4}\left( \frac{1-t}{t}w(t)\right)-4\pi(1-t)^{-1}\right]=-4\pi.\]

Therefore
\[w(t)=4\pi t^2(1-t)^2.\]

Then all the inequalities become equal on $L(0,1)$. That is, for any regular value $t$ of $G$ it follows that the Hessian of $G$ is diagonal in the frame $e_{1}=\frac{\nabla G}{|\nabla G|}, e_{2},e_{3}$, and
\begin{equation}\label{eqn-G11}
  |G_{11}|=|\nabla |\nabla G||.
\end{equation}
Furthermore, by inequality (\ref{eqn-w'}), there exists $\lambda(t)\in \mathbb{R}$ so that
\begin{equation}\label{eqn-nablaG}
  \nabla |\nabla G|=\lambda(t)\nabla G \text{ on } l(t) 
\end{equation}
and $|\nabla G|$ depends only on $t$ on $l(t)$. By inequality (\ref{eqn-h})
\begin{equation}\label{eqn-dw}
  \frac{dw}{dt}(t)=\frac{2(1-2t)}{t(1-t)}w(t)=2h(t)w(t).
\end{equation}
By inequality (\ref{eqn-scalar}), the scalar curvature $\operatorname{Sc}_{g}=0$.
Along the gradient flow $\frac{d\Phi}{dt}=\frac{\nabla G}{|\nabla G|^2}$, we have $\frac{d|\nabla G|}{dt}=\lambda(t)$. From equation (\ref{eqn-dw-w}), we can obtain $w'(t)=\frac{\lambda(t)}{|\nabla G|}w(t)$. Hence
\begin{eqnarray}\label{eqn-GG}
  \frac{d|\nabla G|}{dt}=\frac{w'(t)}{w(t)}|\nabla G|.
\end{eqnarray}
And then $|\nabla G|=cw(t),\; \lambda(t) = cw'(t)$ for some constant $c>0$. Combining Lemma \ref{lem-ineq-u}, we have the following results:
\begin{eqnarray*}
  G_{11}&=&\lambda(t) |\nabla G|,\; G_{12}=G_{13}=0\\
  G_{22}&=&G_{33}=-\frac{1}{2}\lambda(t)|\nabla G|,\; G_{23}=0\\
  (A_{t})_{22}&=&(A_{t})_{33}=-\frac{1}{2}\lambda(t),\; (A_{t})_{23}=0 \text{ on } l(t)\\
  H_{t}&=&-\frac{G_{11}}{|\nabla G|}=-\lambda(t) \text{ on } l(t).
\end{eqnarray*}

Since $|\nabla G|=cw(t)$, we know that $|\nabla G|\neq 0$ for any $t\in (0,1)$. Considering the integral flow of $|\nabla G|^{-2}\nabla G$, the metric $g$ on $(0,1)\times \mathbb{S}^{2}$ can be written as
\[g=\frac{1}{|\nabla G|^2}dt^2+g_{t}=\frac{1}{c^2w^2(t)}dt^2+g_{t},\]
where $g_{t}$ is a metric induced from $g$ on $\Sigma_{t}=\{t\}\times \mathbb{S}^{2}$. The derivative of $g_{t}$ with respect to $t$ satisfies
\begin{eqnarray*}
  \frac{d g_{t}}{d t} &=& \frac{2}{|\nabla G|}A_{t}=\frac{2}{cw(t)}A_{t}=-\frac{\lambda(t)}{cw(t)}g_{t}=-\frac{w'(t)}{w(t)}g_{t}.
\end{eqnarray*}
Then, for some constant $c_1\neq 0$,
\[g_{t}=\frac{1}{c^2_{1}w(t)}g_{\frac{1}{2}}.\]
Therefore
\[g=\frac{1}{c^2w^2(t)}dt^2+\frac{1}{c^2_{1}w(t)}g_{\frac{1}{2}}=\frac{1}{c^2}\left[\frac{1}{w^2(t)}dt^2+\frac{c^2}{c^2_{1}w(t)}g_{\frac{1}{2}}\right].\]
By the scalar flat condition, we can obtain
\[c_{1}^{-2}g_{\frac{1}{2}}=c^{-2}g_{\mathbb{S}^2}.\]
Thus
\[g=\frac{1}{c^2w^2(t)}dt^2+\frac{1}{c^2w(t)}g_{\mathbb{S}^2},\; \text{ on } (0,1)\times \mathbb{S}^{2}.\]

\begin{remark}\label{Rmk-fh}
  In equation (\ref{eqn-f-h}), if we let $f(r)=r^{-2}$. Then equation (\ref{eqn-f-h}) becomes
\[6r^{-4}+3(-2r^{-3}h(r)+r^{-2}h'(r))+3r^{-2}h^2(r)=0.\]
I.e.,
\[r^{2}h'(r)-2rh(r)+r^2h^2(r)+2=0.\]
We can solve this differential equation.
Let $h(r)=\frac{u(r)}{r}$, then
\[h'(r)=\frac{u'(r)r-u(r)}{r^2}.\]
Thus we obtain
\[u'(r)r-u(r)-2u(r)+u^2(r)+2=0.\]
Using the method of separation of variables, we have
\[h(r)=\frac{1-2cr}{r-cr^2}.\]

Munteanu-Wang\cite{WM-IMRN-2023} addressed the case $c=0$, corresponding to $h(r)=r^{-1}$. In the present work, we focus on the case $c=1$, i.e.,
\[h(r)=\frac{1-2r}{r(1-r)}.\]
For Theorem \ref{Thm-c=-1}, we consider the case $c=-1$. The proof of Theorem \ref{Thm-c=-1} follows a similar line of reasoning and is thus omitted herein. A natural question arises: in Equation (\ref{eqn-f-h}), are there other meaningful choices of $f$ and $h$ that yield scalar-flat metrics?
\end{remark}

\bibliographystyle{plain}
\bibliography{harmonic0903.bib}

\end{document}